\magnification=1200
\overfullrule=0pt
\hsize 15.7 truecm \hoffset 0.5 truecm
\vsize 23.5truecm \voffset =0  truecm
\voffset=2\baselineskip
\baselineskip=16pt plus 1pt
\parskip=5pt plus 1pt


\def\BBF#1{\expandafter\edef\csname #1#1\endcsname{%
   {\mathord{\bf #1}}}}
\BBF m
\BBF C
\BBF N
\BBF P
\BBF Q
\BBF R
\BBF T
\BBF Z

\def\bR{{\bf R}}
\def\bZ{{\bf Z}}

\def\cA{{\cal A}}
\def\cE{{\cal E}}
\def\cF{{\cal F}}
\def\cG{{\cal G}}
\def\cIH{{\cal IH}}

\def\CAL#1{\expandafter\edef\csname #1\endcsname{%
   {\mathord{\cal #1}}}}
\CAL A
\CAL C
\CAL E
\CAL F
\CAL G
\CAL H
\CAL{IH}
\CAL K
\CAL P

\def\cA{{\cal A}}
\def\cE{{\cal E}}
\def\cF{{\cal F}}
\def\cG{{\cal G}}
\def\cIH{{\cal IH}}

\def\b{^{\scriptscriptstyle\bullet}}
\def\longto{\longrightarrow}

\def\phi{\varphi}
\def\MOPnl#1{\expandafter\edef\csname #1\endcsname{%
   {\mathop{\rm #1}\nolimits}}}
\MOPnl{odd}
\MOPnl{or}
\MOPnl{Tor}


\vglue 2truecm
\noindent {\bf TOWARDS A COMBINATORIAL INTERSECTION COHOMOLOGY FOR FANS}
\bigskip

\noindent{\bf Karl-Heinz FIESELER} \vskip 6pt

{\parindent=3mm
Matematiska Institutionen,
Box 480,
Uppsala Universitet,
SE-75106 Uppsala, Su\`ede

E-mail: khf@math.uu.se
\par}

\medskip

{\parindent=2cm\rightskip 1.5cm\baselineskip=9pt
\item{\bf Abstract.~}{
The real intersection
cohomology $IH\b(X_\Delta)$ of a
toric variety  $X_\Delta$ is
described in a purely combinatorial way
using methods of
elementary commutative algebra only. We define, for arbitrary fans,
the notion of a ``minimal extension sheaf''
$\cE\b$ on the fan $\Delta$ as an
axiomatic characterization of  the
equivariant intersection cohomology
sheaf. This provides a
purely algebraic interpretation of the
$f$- and $g$-vector of an arbitrary
polytope or fan under a natural
vanishing condition. --- The results
presented in this note originate from
joint work with G.Barthel, J.-P.Brasselet
and L.Kaup (see [3]).\
}} \vskip 10pt

%
%

\noindent{\bf 1. Minimal Extension Sheaves}

\medskip
We introduce the notion of a minimal extension
sheaf on a fan and study some elementary properties of such sheaves.

Let $\Delta$ be a fan in a real vector space $V$ of dimension $n$. We
endow $\Delta$ with the fan-topology, having the subfans $\Lambda$ of
$\Delta$ as the non-empty open subsets. In particular, affine subfans,
i.e. fans $\langle \sigma
\rangle$ consisting of $\sigma$ and its faces, are open;
for simplicity  we usually write $\sigma$ instead of $\langle \sigma
\rangle$. A sheaf $\cF$ of real vector spaces on
$\Delta$ is determined by the collection of  vector spaces $F_\sigma:=\cF
(\sigma)$ for each $\sigma \in \Delta$ together with the restriction
homomorphisms $F_\sigma \to F_\tau$ for $\tau \preceq \sigma$. An
important example is the sheaf $\cA\b$ with $\cA\b
(\sigma):=A\b_\sigma:=S\b(V^*_\sigma)$,
where $V_\sigma:=$span$(\sigma)$, and the natural restriction
homomorphisms. Its sections over a subfan $\Lambda$
are the piecewise polynomial functions on the support $|\Lambda|$ of
$\Lambda$.
\par

\noindent
{\bf Warning:} We use a topologically motivated grading for
$S\b(V^*_\sigma)$: Linear polynomials are of degree $2$, etc.

We set $A\b:=S\b(V^*)$. For a graded $A\b$-module $F\b$, let
$\overline F\b:= R\b \otimes _{A\b} F\b = F\b/ \mm F\b$ denote
the residue class vector space modulo $\mm:=A^{>0}$, where $R\b :=
A\b/\mm$.

\par

\medskip\noindent
{\bf Definition:} {\it A sheaf $\cE\b$ of
graded $\cA\b$-modules on the fan $\Delta$
is called a {\it minimal extension
sheaf\/} (of~$\RR\b$) if it satisfies the
following conditions:

\smallskip{
\def\litem{\par\noindent
                   \hangindent=36pt\ltextindent}
\def\ltextindent#1{\hbox to \hangindent{#1\hss}\ignorespaces}
\litem{(N)} {\it Normalization:\/} One has $E\b_o \cong A\b_o = \RR\b$
for the zero cone~$o$.
\litem{(PF)} {\it Pointwise Freeness:\/} For each cone $\sigma \in
\Delta$, the module~$E\b_{\sigma}$ is {\it free} over~$A\b_{\sigma}$.
\litem{(LME)} {\it Local Minimal Extension $\bmod\ \mm$:\/} For each
cone $\sigma \in \Delta \setminus \{ 0 \}$, the restriction mapping
$
  \phi_\sigma \colon E\b_{\sigma} \to E\b _{\partial \sigma}
$
induces an isomorphism
$
  \overline \phi_\sigma \colon
\overline{E}\b_{\sigma}
  \buildrel \cong \over \longrightarrow
  \overline{E}\b_{\partial\sigma}
$
of graded real vector spaces.\par
}}\par

\smallskip
Condition~(LME) implies that $\cE\b$ is minimal in the set of
all flabby sheaves of graded $\cA\b$-modules satisfying conditions~(N)
and~(PF), whence the name
``minimal extension sheaf''. Moreover, $\cE\b$
vanishes in odd degrees.
For a cone $\sigma \in \Delta$ and
 for a subfan $\Lambda \preceq \Delta$, the $A\b$-modules $E\b_\sigma$
and $E\b_\Lambda$
are finitely generated.

If $\Delta$ is a rational fan for some lattice $N \subset V$ of
maximal rank,
then there is an associated toric variety $X_\Delta$ with the
action of an algebraic torus
$\TT \cong (\CC^*)^n$.
Let $IH\b_\TT(X_\Delta)$ denote the equivariant intersection
cohomology of $X_\Delta$ with real coefficients.
In [1], the following
theorem was proved; it has been the starting point to investigate
minimal extension sheaves.

\medskip \noindent
{\bf Theorem.} {\it Let $\Delta$ be a rational fan.
\item{i)} The assignment
$
\cIH\b_\TT\colon \Lambda
\mapsto IH\b_\TT (X_\Lambda)
$
defines a sheaf on the fan space $\Delta$;
it is a minimal extension sheaf.
\item{ii)} If $\cE\b$ is a minimal
extension sheaf on $\Delta$ and $\sigma $ a $k$-dimensional
cone, then for the local intersection
cohomology $\cIH\b_x$ of $X_\Delta$ in a
point $x$ belonging to the orbit corresponding to $\sigma$, we have:
$
\cIH\b_x ~\cong~ \overline E\b_\sigma\ .
$
\item{iii)} For a complete fan $\Delta$
or an affine fan $\Delta=\langle \sigma
\rangle$ with a cone $\sigma$ of dimension $n$, one has
$
IH\b(X_\Delta) ~\cong~ \overline E\b_\Delta\ .
$}

\smallskip
The vanishing axiom for local intersection cohomology together with
statement ii) in the above theorem yields that, in the case of
a rational fan, the following vanishing condition is satisfied:
\smallskip

\noindent
{\bf Vanishing Condition V($\sigma$):}
{\it For a cone $\sigma$ and a minimal
extension sheaf $\cE\b$ on a fan $\langle \sigma
\rangle$, we have $\overline E^q_{\sigma}~ =~ 0 {\rm\ ~for~ \
} q \ge \dim \sigma\;.
$ } \par

On every fan $\Delta$
there exists a minimal extension sheaf
$\cE \b$. Furthermore,
for any two such sheaves
$\cE\b$, $\cF\b$ on $\Delta$,
each isomorphism $E\b_o
\cong \RR \b \cong F \b_o$
extends (non-canonically) to an isomorphism
$\cE\b \buildrel \cong \over
\to \cF\b$
of graded ${_\Delta\cA}\b$-modules, which
is unique in the case of a simplicial
fan.
Simplicial fans are easily characterized
in terms of minimal extension
sheaves:
The sheaf ${_\Delta\cA}\b$ is a minimal
extension sheaf if  and only if the fan
$\Delta$ is simplicial.

\bigskip
%
%
%

\noindent{\bf 2. Combinatorial equivariant perverse sheaves}

We propose a definition for (combinatorially) ``perverse''
sheaves, here called semisimple sheaves.\par

\medskip
\noindent
{\bf Definition:} {\it A {\it
(combinatorially) semi-simple sheaf\/}
$\cF\b$ on  a fan space $\Delta$ is a
{\it flabby\/} sheaf of graded
$\cA\b$-modules  such that, for each cone
$\sigma \in \Delta$, the
$A\b_{\sigma}$-module
$F\b_{\sigma} $ is {\it
finitely generated and free}.}

For each cone $\tau \in \Delta$,
we construct inductively a ``simple'' sheaf ${}_{\tau}\cE\b$ on
$\Delta$ as follows:
For $\sigma \in \Delta^{\le \dim \tau}:=\{ \sigma \in \Delta; \dim
\sigma \le \dim \tau \}$ we set
$$
  {}_{\tau}E\b_{\sigma}~ :=~ {}_{\tau}\cE\b(\sigma)
  ~:=~ \cases{ A\b_{\tau} & if $\sigma = \tau$, \cr
             0        & otherwise.\cr }
$$
Now, if ${}_{\tau}\cE
\b$ has been defined on
$\Delta^{\le m}$ for some $m \ge
\dim \tau$,   then for each $\sigma \in
\Delta^{m+1}$, we
set $
  {}_{\tau}E\b_{\sigma} :=
  A\b_{\sigma} \otimes_{\bR} \overline
{}_{\tau}\overline E\b_{\partial \sigma}
$
with the restriction map
${}_{\tau}E\b_{\sigma} \to
{}_{\tau}E\b_{\partial\sigma}$
being induced  by some homogeneous
$\bR$-linear section
$s:{}_{\tau}\overline{E}\b_{\partial\sigma}
\longto {}_{\tau}E\b_{\partial \sigma}$ of
the residue class map
${}_{\tau}E\b_{\partial \sigma} \to
{}_{\tau}\overline E\b_{\partial\sigma}$.

\medskip\noindent
{\bf Decomposition Theorem:} {\it Every
semi-simple sheaf $\cF\b$ on  $\Delta$ is
isomorphic to a finite direct sum $
  \cF\b \cong \bigoplus_{i}
{}_{\tau_{i}}\cE\b[-{\ell}_{i}]^{n_i}
$
of shifted simple sheaves with uniquely
determined cones $\tau_{i}
\in \Delta$, natural numbers $n_i \ge 1$
and integers ${\ell}_{i} \in \bZ$.}

From the theorem in section 3 we then obtain the following consequence:

\medskip\noindent
{\bf Corollary 1:} {\it Let
$\pi\colon \hat \Delta \to \Delta$ be a refinement map of fans with
minimal extension sheaves $\hat \E\b$ resp. $\E\b$.
Then the direct image sheaf $\pi_*(\hat \E\b )$ is a semisimple sheaf,
in particular there is a decomposition $
\pi_*(\hat \E\b) \cong \E\b \oplus
\bigoplus_{i}
{}_{\tau_{i}}\cE\b[-{\ell}_{i}]^{n_i}
$
with cones
$\tau_i \in \Delta^{\ge2}$
and positive integers ${\ell}_i,n_i $.}

\medskip
\noindent
{\bf Corollary 2:} {\it For a simplicial refinement $\hat
\Delta$ of $\Delta$,  let $\hat \A\b$ be the
sheaf of $\hat \Delta$-piecewise
polynomial functions on $
\Delta$. Then a minimal extension sheaf
${_\Delta\E}\b$ on $\Delta$ can be realized as a
subsheaf of ${_\Delta\hat \A}\b$.}

\bigskip
%
%

\noindent{\bf 3. Cellular Cech Cohomology of Minimal Extension Sheave}

\medskip

In this section,
we investigate under which assumptions
the module of global sections
$E\b_\Delta := \cE\b(\Delta)$ of a minimal
extension sheaf $\cE\b$ is a free
$A\b$-module.

\medskip
\noindent
{\bf Definition:} {\it A fan $\Delta$ is called {\it
quasi-convex\/} if
for a minimal extension sheaf $\cE\b$ on
$\Delta$, the
$A\b$-module $\cE\b(\Delta)$ is free.}
\par

According to Proposition 6.1 in [1],
a rational fan $\Delta$ is quasi-convex
if and only if the
intersection cohomology of the associated
toric variety $X_\Delta$ vanishes in odd
degrees.
--- The main tool to be used in the sequel is the ``complex of
cellular cochains with coefficients in
$\cE\b$": To a sheaf $\cF$ of real vector
spaces on the fan $\Delta$, we associate its {\it  ``cellular
cochain complex''} $C\b(\Delta, \cF)$. The
cochain module in degree~$k$ is $\bigoplus_{\dim \sigma = n-k} \cF
(\sigma)$, the coboundary operator $\delta^k\colon C^k(\Delta, \cF) \to
C^{k+1}(\Delta, \cF)$ is defined with respect to fixed orientations as
in the usual \v{C}ech cohomology. For $\cF =\cE\b$, the above complex is
-- up to a rearrangement of the indices -- a
``minimal complex'' in the sense of Bernstein and Lunts [4].

We also have to
consider relative cellular cochain
complexes with respect to the boundary
subfan $\partial \Delta$ of a purely $n$-dimensional fan $\Delta$,
supported by the topological boundary of $\Delta$.

\noindent
{\bf Definition:} {\it If the fan
$\Delta$ is purely $n$-dimensional, then
for a sheaf $\cF$ of real
vector spaces on $\Delta$
we set
$$
C\b(\Delta, \partial
\Delta; \cF) :=
C\b(\Delta; \cF)/C\b(\partial \Delta; \cF)
\ ,
$$
where $C\b(\partial \Delta; \cF)
\subset C\b(\Delta; \cF)$ is the
subcomplex of cochains supported in
$\partial \Delta$.}

We also need the augmented
complex
$$
\tilde C\b(\Delta, \partial
\Delta; \cF):
0 \longto \cF(\Delta)
\longto C^0(\Delta, \partial
\Delta; \cF) \longto \dots  \longto
C^n(\Delta, \partial
\Delta; \cF) \longto 0
$$
and its cohomology groups
$
\tilde H^q(\Delta, \partial
\Delta; \cF) :=
H^q(\tilde C\b(\Delta, \partial
\Delta; \cF))\ .
$
\medskip

\noindent
{\bf Theorem:} {\it For a purely
$n$-dimensional fan $\Delta$ and a minimal extension
sheaf~$\cE\b$ on~$\Delta$, the following statements are equivalent:
\item{i)} We have
$
\tilde H\b(\Delta, \partial \Delta;
\cE\b)= 0$.
\item{ii)} The $A\b$-module $E\b_\Delta:=\cE\b(\Delta)$ of
global sections is free.
\item{iii)} The support $|\partial
\Delta|$ of the boundary subfan is a real homology manifold.
\par
\noindent
For a rational fan $\Delta$,
the above conditions are equivalent to
\item{iv)}
For the toric variety
$X_\Delta$ associated to $\Delta$, we have
$IH^{\odd}(X_\Delta)=0$.}

\smallskip

Since complete fans are quasi-convex, the previous results
provide a proof of a conjecture of
Bernstein and Lunts (see [4],~15.9).

\noindent
{\bf Corollary :} {\it For a complete fan $\Delta$,
the minimal complex of Bernstein and Lunts is exact.}

Furthermore, for a quasi-convex fan~$\Delta$ and
a minimal extension sheaf $\cE\b$ on~$\Delta$,
even the $A\b$-submodule $E\b_{(\Delta, \partial \Delta)}$ of
$E\b_\Delta$  consisting of the global sections vanishing on the
boundary subfan $\partial \Delta$ is a free $A\b$-submodule.

\bigskip
%
%

\noindent{\bf 4. Poincar\'e Polynomials and Poincar\'e duality}

\medskip

For a quasi-convex fan~$\Delta$ and a a minimal extension sheaf $\cE\b$
on~$\Delta$, we want to discuss the Poincar\'e polynomials related to
$\Delta$ and the pair $(\Delta, \partial\Delta)$.

\noindent
{\bf Definition:} {\it The Poincar\'e polynomial of
$\Delta$ is the polynomial $
  P_{\Delta} (t) := \sum_{q \ge 0}^{<\infty} \dim\,
 \overline E_{\Delta}^{2q}
  \cdot t^{2q}$.

\noindent
The relative Poincar\'e polynomial $P_{(\Delta, \partial \Delta)}$
is defined in an analogous manner.
}
\par

The relation between a global Poincar\'e polynomial $P_\Delta$ and its
local Poincar\'e polynomials $P_\sigma$ for $\sigma \in \Delta$ is rather
explicit:

\medskip
\noindent
{\bf Local-to-Global
Formula:} {\it If $\Delta$ is a
quasi-convex fan of dimension~$n$,
we have
$$
  P_{\Delta}(t) = \sum_{\sigma \in
\Delta \setminus \partial \Delta}
  (t^2-1)^{n- \dim \sigma} P_{\sigma}(t)
  \quad\hbox{and}\quad
P_{(\Delta, \partial \Delta)}(t)
=\sum_{\sigma \in
\Delta}
  (t^2-1)^{n- \dim \sigma} P_{\sigma}(t)
  \ .
$$}

\medskip
The proof of the above formul{\ae} depends
on the fact that $C\b(\Delta, \partial
\Delta; \cE\b)$ resp. $C\b(\Delta;\cE\b)$
are resolutions of $E\b_\Delta$ resp.
$E\b_{(\Delta, \partial \Delta)}$. Hence,
the Poincar\'e series of  $E\b_\Delta$ resp.
$E\b_{(\Delta, \partial \Delta)}$ equals the alternating sum of the
Poincar\'e series of the cochain
modules $C^i(\dots)$. Finally, we use the fact
that $E\b_\Delta \cong A\b \otimes_\RR
\overline E\b_\Delta$, since
$E\b_\Delta$ is free, and similarly for
$E\b_{(\Delta, \partial \Delta)}$, while
$E\b_\sigma \cong A\b_\sigma \otimes_\RR
\overline E\b_\sigma$.
\par
Using an  induction argument one proves:
\medskip
\noindent
{\bf Corollary:} \hangafter=1\hang {\it
Let $\Delta$ be a quasi-convex fan.
\item{i)} The relative Poincar\'e polynomial
$P_{(\Delta, \partial \Delta)}$ is monic of degree~$2n$.
\item{ii)} The absolute Poincar\'e polynomial
$P_\Delta$ is of degree $2n$ iff $\Delta$ is complete; otherwise, it
is of strictly smaller degree.
\item{iii)} For a non-zero cone $\sigma$, the local Poincar\'e polynomial
$P_\sigma$ is of degree at most $2 \dim \sigma~-~2$.}

Of course, statement ii) is a rather weak vanishing
estimate; in fact we expect the much stronger
vanishing condition V($\sigma$) to hold.
\par

In order to have a recursive
algorithm for the computation of global
Poincar\'e polynomials, we relate
the local Poincar\'e polynomial $P_\sigma
$ to the global one
of some fan
$\Lambda_\sigma$ in a vector space of
lower dimension: The fan $\Lambda_\sigma$ ``lives''
in the quotient vector space
$V_\sigma/L$, where~$L$ is a line
 in~$V$ passing through the relative
interior of~$\sigma$.
For the projection $\pi\colon V_\sigma \to V_\sigma/L$, we pose $
\Lambda_\sigma := \{ \pi(\tau); \tau
\prec \sigma \}
$.
The homeomorphism
$
\pi|_{\partial \sigma}\colon |\partial \sigma|
\longrightarrow V_\sigma/L
$
induces an
isomorphism of the fans $\partial \sigma$ and
$\Lambda_\sigma$. Using the truncation
operator $\tau_{<j}(\sum a_q t^q)=
\sum_{q <r} a_qt^q$ we can now formulate the
next step:

\medskip
\noindent
{\bf Local Recursion Formula:} {\it Let $\sigma$ be a non-zero cone.
\item{i)} If $\sigma$ is simplicial, then we have $P_\sigma \equiv 1$.
\item{ii)} If the condition $V(\sigma)$ of  section 1 is satisfied, then
we have
$
  P_{\sigma}(t) =
  \tau_{<\dim\sigma}\big((1-t^2)
  P_{\Lambda_\sigma}(t)\big)$.}

For the proof, we consider a minimal
extension sheaf $\cG\b$ on the fan
$\Lambda:=\Lambda_\sigma$. Let $B\b$
be the polynomial algebra on $V_\sigma /L$, considered as subalgebra of
$A\b_\sigma \cong B\b[T]$.
Then we have an $B\b$-module isomorphism
$G\b_\Lambda \cong E\b_{\partial
\sigma}$, such that $\overline
E\b_{\partial \sigma} \cong \overline
G\b_\Lambda / T \overline
G\b_\Lambda$: Here $\overline
G\b_\Lambda$ is the residue class module
of the $B\b$-module $G\b_\Lambda$ and $T$
acts on it via the isomorphism
$G\b_\Lambda \cong E\b_{\partial
\sigma}$, the latter module living over
$A\b_\sigma$. The action of $T$ on $G\b_\Lambda$
coincides with the multiplication with
the piecewise linear strictly convex
function $\psi:= T \circ (\pi|_{\partial
\sigma})^{-1} \in \cA^2 (\Lambda)$. Now
we use the following combinatorial
version of the Hard Lefschetz Theorem:
\medskip

\noindent
{\bf Combinatorial Hard Lefschetz Theorem:} {\it In the
same notations  as
in the proof of the Local Recursion Formula we set
$m:= \dim (V_\sigma /L)=\dim \sigma -1$. If
the condition $V\big(\sigma)$ is
satisfied, then $ \mu\colon G\b_{\Lambda}
\to G\b_{\Lambda}[2] \ , \ f \mapsto \psi
f $ induces a  map $ \overline \mu^{2q}:\overline
G^{2q}_{\Lambda} \longrightarrow \overline
G^{2q+2}_{\Lambda}\;, $
which is injective for $2q \le m-1$ and
surjective for $2q \ge m-1$.}
\medskip
The surjectivity is nothing but a
reformulation of the vanishing condition
V($\sigma$), whence the injectivity is
obtained via Poincar\'e duality for
the real vector space $\overline
G\b_\Lambda$, the map $\overline \mu$
being selfadjoint with respect to the
Poincar\'e duality pairing:

Such a Poincar\'e duality on a quasi-complex fan $\Delta$
is obtained as follows: By a stepwise procedure,
one constructs an internal (non-canonical) intersection product
$\cE\b \times \cE\b \to \cE\b$. Then one composes
the induced product on the level of global sections with an
evaluation map $E\b_{(\Delta, \partial \Delta)} \to A\b[-2n]$,
which is homogeneous of degree 0 and unique up to a non-zero real
factor. That construction uses the above corollary and the freeness of
$E\b_{(\Delta, \partial \Delta)}$. The pairing thus obtained
induces a pairing on the level of residue class vector spaces.

\medskip
\noindent
{\bf  Poincar\'e Duality Theorem:} {\it For every
quasi-convex fan $\Delta$, the pairings
$$
\eqalign{
E\b_\Delta
\times
E\b_{(\Delta, \partial \Delta)}
&\longrightarrow
E\b_{(\Delta, \partial \Delta)}
\longrightarrow
A\b[-2n] \cr
\overline E\b_\Delta
\times
\overline E\b_{(\Delta, \partial \Delta)}
&\longrightarrow
\overline E\b_{(\Delta, \partial \Delta)}
\longrightarrow
\RR\b[-2n]\ .\cr}\leqno{\it and}
$$
are dual pairings of free $A\b$-modules
resp. of $\RR$-vector spaces.
}

We end this section with
a numerical version of Poincar\'e duality:

\noindent
{\bf Corollary:} {\it For a quasi-convex fan~$\Delta$,
the global Poincar\'e polynomials $P_\Delta$ and
$P_{(\Delta,\partial\Delta)}$ are related by the
identity
$$
P_{(\Delta,\partial \Delta)}(t)~=~t^{2n} P_\Delta(t^{-1})\; .
$$}


\centerline{\bf R\'ef\'erences bibliographiques}

\medskip

{
\item{[1]} Barthel~G., Brasselet~J.-P., Fieseler~K.-H., Kaup~L.,
Equivariant Intersection Cohomology of Toric Varieties,
Algebraic Geometry: Hirzebruch 70. Contemp. Math. AMS~241,
(1999) 45--68.
\smallskip

\item{[2]} ---,
Equivariant Intersection
Cohomology of Toric Varieties,
Dep. Math. Uppsala, U.U.D.M. Report 1998:34.
\smallskip

\item{[3]} ---, Combinatorial Intersection Cohomology of Fans,
To appear.
\smallskip

\item{[4]} Bernstein~J., Lunts~V., Equivariant
Sheaves and Functors, Lecture Notes in Math., vol. 1578,
Springer-Verlag
Berlin etc., 1993.
\smallskip

\item{[5]}  Brion~M., The Structure of the Polytope
Algebra, T\^ohoku Math. J., 49, (1997), 1--32.
\smallskip

\item{[6]} Brylinski~J.-L., Equivariant Intersection
Cohomology, in: Kazhdan-Lusztig Theory and Related Topics,
Contemp.\ Math.\ vol. 139, Amer.\ Math.\ Soc., Providence, R.I.,
(1992), 5--32.
\smallskip

\item{[7]} Fieseler~K.-H., Rational Intersection
Cohomology of Projective Toric Varieties, Journ. reine angew. Math.
(Crelle), 413, (1991), 88--98.
\smallskip

\item{[8]} Goresky~M., Kottwitz~R., MacPherson~R.,
Equivariant Cohomology, Koszul Duality, and the Localization Theorem,
Invent. Math., 131, (1998), 25--83.
\smallskip

\item{[9]} Oda~T., The Intersection Cohomology and
Toric Varieties, in:  T.~Hibi (ed.), Modern Aspects of
Combinatorial Structure on Convex Polytopes, RIMS Kokyuroku,
857 (Jan. 1994), 99--112.
\smallskip
}
\end